\documentclass[aps,12pt,nofootinbib]{revtex4}

\usepackage{amsbsy}

\usepackage{graphicx,graphics}
\usepackage{enumerate}
\usepackage{subfigure}
\usepackage{amsmath}
\usepackage{amsfonts}
\usepackage{amssymb}
\usepackage{amsthm}
\usepackage{psfrag}

\catcode`\@=11
\def\numberbysection{\@addtoreset{equation}{section}
\def\theequation{\arabic{section}.\arabic{equation}}}
\numberbysection

\def\Z{\mathbb{Z}}

\begin{document}

\title{Perfect boundaries in rotor-router aggregation on cylinders}

\author{V.B. Priezzhev}
\affiliation{Bogoliubov Laboratory of Theoretical Physics,\\ Joint Institute for Nuclear Research, 141980 Dubna, Russian Federation}

\begin{abstract}
We study a rotor-router version of the internal diffusion-limited aggregation introduced by J.Propp. The existing estimations of boundary fluctuations of the aggregation cluster show that they grow not faster than $O(\log r)$ with the cluster radius $r$. We consider the rotor-router internal DLA on the semi-infinite cylinder and prove a constant width of boundary fluctuations not depending on the radius of the cylinder.
\end{abstract}

\maketitle


\section{Introduction}
The internal diffusion-limited aggregation (internal DLA) was introduced by Meakin and Deutch \cite{Meakin} as a model of cluster growth phenomena.
Given a connected cluster of $N$ sites $A(N)\subset \Z^d$ containing the origin, let $A(N+1)$ be a cluster obtained by a random walk starting at the origin and stopping when the walker first exits $A(N)$, so that the site $A(N+1)\backslash A(N)$ is the endpoint of the walk. After each stop, the walker returns to the origin and the process continues. If one thinks of the random walk as a trail of a diffusing particle, $A(N)$ is the cluster of $N$ particles delivering from an unrestricted source at the origin. The question is the shape of the cluster $A(N)$ and the smoothness of its boundary for large $N$ if one starts from $N=1$.

Lawler et al \cite{Lawleretal} showed that the rescaled shape of $A(N)$ is a ball for any dimension. In \cite{Lawler}, Lawler proved that fluctuations of boundaries grow with the ball radius $r$ not faster than $C r^{1/3}$ where $C$ is a constant. Asselah and Gaudilli\`{e}re \cite{AG1,AG2} and independently Jerison,Levine and Sheffield \cite{JLS1,JLS2} presented the new bound $C \log r$ on the maximum fluctuations for two dimensions and $C_d\surd\overline{\log r}$ for $d\geq 3$.

It is intuitively clear that the smoothness of boundaries in the internal DLA depends on statistical properties of the walk inside the cluster of the visited sites. Propp was the first who proposed to replace the simple random walk in the two-dimensional aggregation by a deterministic process called the rotor-router walk \cite{Kleber}. The rotor mechanism introduced in \cite{PDDK} works as follows. At each site of the square lattice there is a rotor pointing to one of four cardinal directions N,E,S,W. At each time step, the walker arriving at a lattice site rotates the rotor in this site clockwise by 90 degrees and takes a step in the new direction. The rotor rule leads to essential changes of properties of the walk trajectories, and as a consequence, to changes in the aggregation process.

The rotor-router internal aggregation for the $d$-dimensional lattice was considered by Levine and Peres, who proved that the asymptotic shape of the cluster remains a ball as in the internal DLA \cite{LP05}.  The in-radius of the set of sites occupied by the aggregation cluster is at least $r-O(\log r)$ and the
out-radius is at most $r+O(r^\alpha)$ for any $\alpha > 1-1/d$ \cite{LP07}. Recently, Levine and Peres \cite{LP16} have presented the stronger result $r+O(\log r)$ for the outer bound as well. These bounds can be compared with the results for the clusters in the divisible sandpile
formed from a pile of a given mass at the origin. Fluctuations of the radius $r$ of the domain occupied by sand are restricted by constants that do not depend on $r$. In \cite{LP07}, Levine and Peres asked the question whether the similar estimations hold for the rotor-router aggregation where the near perfect circularity was observed in numerical simulations. In this paper, we will prove that the boundary fluctuations of a constant width take place in a rotor-router aggregation model for a specific geometry, namely, for a semi-infinite cylinder.

The standard internal DLA for cylinders was analysed by Jerison,Levine and Sheffield in \cite{JLS3}. The cylinder was taken as the square $L\times\infty$ lattice that is periodic in the horizontal direction and extends to infinity in both vertical directions. In the initial state, the half of cylinder with the vertical coordinates $y<0$ is occupied by particles. New particles arrive from a source at $y=-\infty$. The motion of each particle is a random walk on the occupied region finishing in the first unoccupied site. When $N$ particles are added, the occupied region reaches average level $y=N/L$. The authors proved that the maximum deviation of the boundary of the occupied region from the average level is confined in a band of width $C\log L$ with probability tending to 1 when $L\rightarrow \infty$.

We consider a rotor-router version of the internal DLA for the cylinder with slightly different boundary conditions. Instead of the infinite cylinder, we consider a semi-infinite cylinder with $y\geq 0$. All $L$ boundary sites at the level $y=0$ are connected with a single site $S$ by the edges directed from $S$ to the boundary. Four edges emanated from each bulk site and three edges from each boundary site are ordered cyclically, and $L$ edges from $S$ are equipped by an arbitrary cyclic order. Particles one by one start at $S$, jump to one of the $L$ bottom sites and perform the rotor-router walk on the set of previously occupied sites until reaching a first empty site. Then a new particle is added to vertex $S$ which sends it along the next edge in the fixed permutation of $L$ edges directed to the boundary.

The replacement of the random walk by the rotor-router one gives a drastic smoothing of the upper boundary of the region $A(N)$ occupied by $N$ particles. We will show that for each $N=0(mod L)$ the boundary is a perfect ring i.e. the section of the cylinder by a horizontal plane at the level $y=N/L$.
For $N=n(mod L)$, $n$ particles occupy $n<L$ sites on the level $y=\lfloor N/L \rfloor +1$, where specific positions of $n$ sites depend on the random initial orientations of rotors on the cylinder. Thus, the width of the random boundary for $N \neq 0(mod L)$ is exactly 1.

The perfect boundary of $A(N)$ for particular values of $N$ is not the first example of vanishing fluctuations in the aggregation models.
Kager and Levine \cite{Kager} proposed a model of rotor-router aggregation on the square lattice obtained from $\Z^2$ by reflecting all directed edges on the $x$-and $y$-axes that point to a vertex closer to the origin. They proved that there exists a rotor configuration $\rho_0$ such that the occupied domain $A(N)$ for $N=2n(n+1)$ and integer $n$ is exactly the diamond i.e. the set $\{(x,y)\in \Z^2:|x|+|y|\leq n \}$. The difference of our model from the example in \cite{Kager} is that the model considered below gives the perfect boundaries for an arbitrary initial configuration of rotors $\rho_0$ and there is no layered structure of rotor configurations in the domain $A(N)$.

In Section 2, we define the model of rotor-router aggregation on the semi-infinite cylinder and prove the main result. In Section 3, we consider the structure of recurrent states, first for the unicycles and then for arbitrary configurations of rotors on finite sublattices of the semi-infinite cylinder. Section 4 contains a discussion.

\section{The model and main results.}

Consider a part of the square lattice $G\subset \Z^2$ whose vertices are denoted by the pair of integers $(x,y)$, $x=1,2,\dots,L$ with $y\geq 0$ and neighboring vertices are connected by directed edges.
The graph $G$ is periodic in the horizontal direction, so that vertices $(1,y)$ and $(L,y)$ are neighbors for all $y$. We add to $G$ an extra vertex $S$ and connect it with boundary vertices $(1,0),(2,0),\dots,(L,0)$ by edges directed from $S$ to the boundary. To define the the rotor-router dynamics, we fix the clockwise ordering of edges N,E,S,W directed from bulk vertices $(x,y)$, $y>0$, the ordering W,N,S of the edges directed from boundary vertices and an arbitrary ordering of $L$ edges from vertex $S$ to boundary vertices $(\pi(1),0),(\pi(2),0),\dots,(\pi(L),0)$, where $\pi(i),i=1,2,\dots,L$ is the permutation of $L$ elements $\pi \in S_L$.

A state $(w,\rho)$ of the model is characterized by the position $w\in G\cup S$ of the particle called usually {\it chip} and a configuration of rotors $\rho(v)$ at all vertices $v\in G\cup S$. Each step of the rotor-router walk consists of updating the rotor configuration $\rho(v) \rightarrow \rho(v)^{+}$ and moving the chip from $w$ to the nearest neighbor of $w$ in direction $\rho(w)^{+}$. The operation $\rho(v) \rightarrow \rho(v)^{+}$ is the rotation of the rotor at $w$ to the next edge in the fixed order of edges emanating from $w$.

The initial state of the model is the empty set of occupied vertices $A(0)$ and an arbitrary configuration of rotors $\rho_0$ at bulk vertices.
The rotors at boundary vertices can be directed randomly to W,N or E, and the single rotor at vertex $S$ is directed to one of the boundary vertices chosen randomly. The initial position of a chip is vertex $S$.

The aggregation process consists of an unrestricted sequence of rotor-router walks starting at $S$ and finishing at the first vertex $v^{*}\in G$ not visited
before. The aggregation cluster $A(N)\subset G$ becomes $A(N+1)=A(N)\cup \{v^{*}\}$ after the $(N+1)$-th walk. It is convenient to divide the whole sequence of walks into {\it periods} consisting of $L$ walks and investigate configurations of aggregated particles after each period.

The result of the first period is trivial since each of boundary vertices $(1,0)$, $(2,0), \dots, (L,0)$ receives one particle from $S$ directly. The result of the next period is also simple because the vertices of the second row $(1,1),(2,1),\dots,(L,1)$ receive also one particle each during the second period. This follows from the fact that none of the edges connecting the vertices of the first and the second row can conduct more than one particle during the second period. Indeed, assume that the boundary vertex $(i,0), 1\leq i \leq L$ is the first vertex of the first row which transfers two particles to the vertex $(i,1)$. The rotor at $(i,0)$ should make for this more than three rotations, for instance W$\rightarrow$ N, N$\rightarrow$ E, E$\rightarrow$ W, W$\rightarrow$ N.
Since each rotation at a vertex $v$ follows a particle hit into $v$, the vertex $(i,0)$ should receive al least four particles from its neighbors. One particle goes from $S$, but vertices $(i-1,0)$ and $(i+1,0)$ can deliver not more than one particle each since the vertex $(i,0)$ is the first vertex performing more than three rotations by the assumption. We come to a contradiction and, therefore, $L$ particles emanated from $S$ during the second period are distributed uniformly among vertices of the second row.

We obtained the perfect boundary of the aggregation cluster $A(N)$ for $N=L$ and $N=2L$. To prove the same for any $N=0(mod L)$, we need an induction on the number of periods with the following induction hypothesis.

{\it Induction hypothesis}. Let $G_K\subset G$ be a part of graph $G$ consisting of $K$ rows $(1,y),(2,y),\dots,(L,y)$, $0 \leq y \leq K-1$.
Assume that there exists the aggregation cluster $A(N)$ with $N=KL$ obtained after $K$ periods and coinciding with $G_K$. Then, the $(K+1)$-th period results in the uniform distribution of $L$ particles among the vertices of the $(K+1)$-th row.

{\it Proof}. Chips come in $G_K$ from $S$ one by one and walk from the first row to the $K$-th row where they leave $G_K$ along one of the edges connecting the $K$-th and $(K+1)$-th rows. Assume that two chips used the same edge $e_{K-1,K}(i)$ connecting vertices $(i,K-1)$ and $(i,K)$ during one period. Let $t$ be a moment of time within the $(K+1)$-th period when the second chip passes $e_{K-1,K}(i)$. This implies one of two possibilities.

(a)There is a moment of time $t_1 < t$ of the $(K+1)$-th period when one of the rotors of the $K$-the row, say the rotor at $(i,K-1)$, performs $r=5$ rotations for the first time. Five rotations need five hits of the walk into $(i,K-1)$. Since the rotors at $(j,K-1)$,$j\neq i$ performed less than 5 rotations to the moment $t_1$, the vertex $(i,K-1)$ is visited from its neighbors $(i+1,K-1)$ and $(i-1,K-1)$ not more than once from each side. Therefore, it was visited to the moment $t_1$ not less than $\Delta=3$ times from $(i,K-2)$.

(b) Vertex $(i,K-1)$ is visited from $(i,K-2)$ less than 3 times to the moment $t_1$. Then, it should be visited from its neighbors $(i+1,K-1)$ and $(i-1,K-1)$ more than twice to perform five rotations. One of the vertices $(i+1,K-1)$ and $(i-1,K-1)$ or both should be visited to this end not less than five times taking additional hits either from next neighbors in the $K$-th row or from the $(K-1)$-th row. Since the number of vertices in the row is restricted, there is a vertex $0 \leq j \leq L, j\neq i$ which is visited $\Delta=3$ times from $(j,K-2)$ for the first time before $t_1$.

Therefore, there exist a moment of time $t_1^{'}$, $t_1^{'} <t_1 < t$, within the $(K+1)$-th period and a vertex $(k,K-2)$, $0 \leq k \leq L$ of the $(K-1)$-th row where the rotor sends the chip to the $K$-th row exactly $\Delta=3$ times for the first time.

Consider now rotors in the $(K-1)$-th row at vertices $(1,K-2),(2,K-2),\dots,(L,K-2)$. If one of these rotors sends the chip to the $K$-th row three times during the $(K+1)$-th period, there are two possibilities again.

(a) There is a moment of time $t_2 < t_1$ when the rotor at a vertex $(i,K-2)$,$0 \leq i \leq L$ performs $r=9$ rotations for the first time.  Since the rotors at $(j,K-2)$,$j\neq i$ performed less than 9 rotations to the moment $t_2$, the vertex $(i,K-2)$ is visited from its neighbors $(i+1,K-2)$ and $(i-1,K-2)$ not more than two times from each side. By condition, vertex $(i,K-1)$ performs less than 5 rotations to the moment $t_2$. So, it sends not more than one chip towards $(i,K-2)$. Therefore, vertex $(i,K-2)$ is visited to the moment $t_2$ not less than $\Delta=4$ times from $(i,K-3)$ to provide total 9 visits.

(b) Vertex $(i,K-2)$ is visited from $(i,K-3)$ less than 4 times and not more than once from $(i,K-1)$. Then, it must take additional visits from one of the rotors at $(j,K-2), j\neq i$ in the $(K-1)$-th row. At least one of these rotors has to perform 9 rotations before $t_2$ and there exists a rotor which takes $\Delta=4$ visits from the row $(K-2)$.

Therefore, there exists a moment of time  $t_2^{'} <t_2 $, within the $(K+1)$-th period and a vertex $(k,K-3)$, $0 \leq k \leq L$ of the $(K-2)$-th row where the rotor sends the chip to the $(K-1)$-th row exactly $\Delta=4$ times for the first time.

Continuing the reasoning row by row, we come to a conclusion that one of the rotors of the $(K-n)$-th row at vertices $(1,K-n-1),\dots,(L,K-n-1)$ must perform
$r=4(n+1)+1$ rotations. If it happens at the moment $t_{n+1}$, there exists a moment $t_{n+1}^{'} < t_{n+1}$ and a vertex $(k,K-n-2)$, $0 \leq k \leq L$ where the rotor sends the chip to the $(K-n)$-th row exactly $(n+3)$ times for the first time and performs not less than $r=4(n+2)+1$ rotations.

Now we put $n=K-1$ and consider the first row of rotors at vertices $(1,0),\dots,(L,0)$. It follows from the above that at least one of these rotors should be directed to the second row $K+1$ times.The rotors at boundary of the cylinder have three allowed directions W,N,S so that a full rotation consists of 3 steps. The minimal number of rotations needed to send chips $K+1$ times to the second row is $r=3K+1$. If the rotor at vertex $(k,0),0\leq k \leq L$, rotates the $(3K+1)$-th time at the moment $t_K$ for the first time,it implies that this vertex took $3K+1$ visits from neighboring vertices . Two vertices $(k-1,0)$ and $(k+1)$ can provide $(k,0)$ with no more than $K$ visits each to moment $t_K$, and vertex $(k,1)$ not more than $K-1$ visits, altogether $3K-1$ visits. Therefore, the vertex $S$ should send to $(k,0)$ not less than 2 chips what contradicts the condition that $S$ distributes during one period exactly $L$ chips among all $L$ outgoing edges.

Thus, the induction hypothesis is proved and we obtain the perfect boundary of the aggregation cluster $A(N)$ for any $N=0(mod L)$. If the number of rotor-router walks starting from vertex $S$ has the form $N=n(mod L)$, the aggregation cluster consists of $\lfloor N/L \rfloor$ rows of vertices completely
occupied by particles, and the remaining $n$ particles are arranged in the $(\lfloor N/L \rfloor+1)$-row at $n$ vertices whose positions depend on a random configuration of rotors $\rho(v)$ on the occupied part of the semi-infinite cylinder. The width of the random boundary for $N \neq 0(mod L)$ is 1 as claimed in the Introduction.

\section{Multi-Eulerian tours and multicycles}

In this section, we describe the rotor-router walk during one period in more detail. The aggregation cluster $A(N)$ obtained after $K$ periods consists of
$K$ rows $(1,y),(2,y),\dots,(L,y)$, $0 \leq y \leq K-1$, so that $N=KL$. The part of graph $G$ whose set of vertices $V$ coincides with $A(N)$ is $G_K$.
To consider the $(K+1)$-th period, we construct an auxiliary graph $\overline{G_K }$ in the following way: a) add to $G_K$ an additional vertex $S$ connected with boundary vertices $(1,0),(2,0),\dots,(L,0)$ by $L$ edges directed from $S$ to the boundary; b)add to $G_K$ the second vertex $S^{'}$ connected with the top row $(1,K-1),\dots,(L,K-1)$ by $L$ edges directed towards $S^{'}$; c)join vertices $S$ and $S^{'}$. An example of $\overline{G_K }$ for $L=4$ and $K=3$ is shown in Fig.\ref{uni}. The obtained graph $\overline{G_K }$ is a finite directed non-Eulerian graph. Rotors have three possible directions at vertices of the first row $(1,0),(2,0),\dots,(L,0)$, four directions at vertices of rows
$(1,y),(2,y),\dots,(L,y)$, $1 \leq y \leq K-1$ and $L$ directions at $S$.

We put a single chip to vertex $S$ and consider the path traced by the chip performing the rotor walk until the moment when the chip has used all $L$ outgoing edges of $S$ and returned to $S$. By construction, the closed path on $\overline{G_K }$ is equivalent to $L$ paths traced by $L$ particles starting at $S$ during the $(K+1)$-th period. In the original cylinder formulation, each of these $L$ particles reaches the $(K+1)$-th row and then the process returns to $S$. In the case of $\overline{G_K }$, all vertices of the $(K+1)$-th row are gathered into the single vertex $S^{'}$ which coincides with $S$.

For the closed path on $\overline{G_K }$, a question arises what is the number of passages of each edge in the path. A key role in answering this question is played by the notion of
{\it unicycle} which is defined in \cite{HLMPPW} as the chip-and-rotor state $(w,\rho)$ where a configuration of rotors $\rho(v)$ contains a unique directed cycle and vertex $w$ lies on this cycle. In the case of Eulerian graphs, where the indegree of each vertex equals its outdegree, there exists a tour that traverses each directed edge exactly once \cite{PDDK,HLMPPW}. The condition that the chip performs the Eulerian tour is formulated
as Lemma 4.9 in \cite{HLMPPW}: let $G$ be an Eulerian directed graph with $m$ edges. Let $(w,\rho)$ be a unicycle in $G$. If we iterate the rotor-router operation $m$ times starting from $(w,\rho)$, the chip traverses an Eulerian tour of $G$, each rotor makes one full turn, and the state of the system returns to $(w,\rho)$.

In the case of non-Eulerian graphs, Farrell and Levine \cite{FarrLev} introduced a notion of $\pi$-Eulerian tour to deal with a situation when each edge can be traversed many times.
Given a directed graph $G(V,E)$ with the vertex set $V$ and the edge set $E$, let $\boldsymbol{\pi}=(\pi_1,\pi_2,\dots,\pi_{|V|})$ be a vector with strictly positive integer values.
A $\pi$-Eulerian tour of $G$ is a closed path that uses each directed edge of $G$ exactly $\pi_{i(e)}$ times where $i(e)$ is the tail of edge $e\in E$ \cite{FarrLev}. The theorem proved in \cite{FarrLev} reads that $G$ has $\pi$-Eulerian tours if and only if
\begin{equation}
\Delta \boldsymbol{\pi}=0
\label{Laplace}
\end{equation}
where $\Delta$ is the Laplacian of $G$.

The notion of $\pi$-Eulerian tours allows us to give the description of the rotor-router walk during one period.
It is convenient to divide the rest of this section into two subsections considering separately the case when the chip-and-rotor state $(w,\rho)$ on the graph $\overline{G_K }$ is an unicycle,and the case when the rotors are allowed to form a number of cycles and should be described in terms of multicycles \cite{PPP15}.

\subsection{Unicycles}

\begin{figure}[!ht]
\includegraphics[width=60mm]{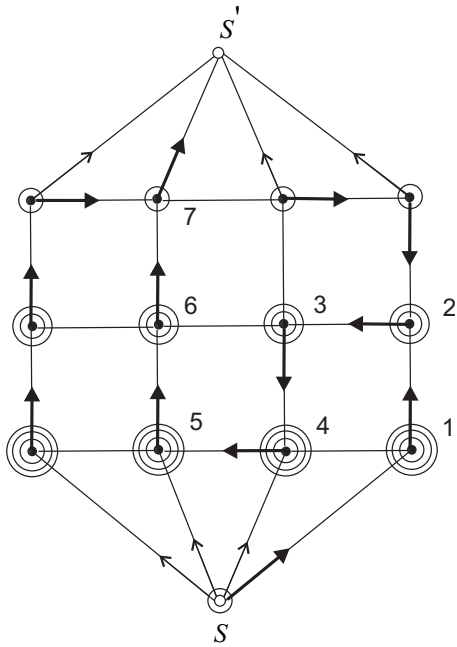}
\vspace{-1mm}
\caption{\label{uni} The auxiliary graph  $\overline{G_K }$ for $K=3,L=4$. Bold vectors denote the initial rotor configuration $\rho_0$. The initial position of chip is $S$. The chip-and-rotor state $(S,\rho_0)$ is a unicycle with cycle $S,1,2,\dots,7,S^{'}$. The rings around the vertices denote multiple full rotations of rotors during the $\pi$-Eulerian tour. }
\end{figure}

An example of a unicycle on the graph $\overline{G_K }$ is given in Fig.\ref{uni}. To relate the unicycle of the non-Eulerian graph to the $\pi$-Eulerian tours, we use a trick outlined in \cite{FarrLev}. Define a multi-graph $\widetilde{G_K}$ replacing each edge  $e$ of $\overline{G_K }$ by $\pi_{i(e)}$  edges $e^1,\dots,e^{\pi_{i(e)}}$. Each vertex $v$ of $\widetilde{G_K}$ has outdegree $d_v \pi_v$ and indegree $\sum_{u\in V} \pi_u d_{uv}$ where $d_v$ is outdegree of $v$ for the graph $\overline{G_K }$ and $d_{uv}$ is the number of edges directed from $u$ to $v$. Then, the graph $\widetilde{G_K}$ is Eulerian if and only if (\ref{Laplace}) holds. If $(e_1^{i_1},\dots,e_m^{i_m})$ is an Eulerian tour of $\widetilde{G_K}$, then $(e_1,\dots,e_m)$
is a $\pi$-Eulerian tour of $\overline{G_K }$. Therefore, we get from Lemma 4.9,\cite{HLMPPW} a conclusion that the chip starting from the initial rotor configuration $\rho_0$ on $\overline{G_K }$ traverses a $\pi$-Eulerian tour if and only if the state $(w,\rho_0)$ is a unicycle. The explicit form of vector $\boldsymbol{\pi}=(\pi_1,\pi_2,\dots,\pi_{|V|})$ depends on the geometry of graph $\overline{G_K }$ and can be easily found by use of the results of the previous section. Indeed, the absence of repeated exits from the vertices of the row $(1,K-1),(2,K-1),\dots,(L,K-1)$ along $L$ edges directed to $S^{'}$ suggests the choice $\pi_{i_{K-1}}=1$ for all vertices of the top row $(i_{K-1},K-1), i_{K-1}=1,\dots,L$. The choice $\pi_S=1$ follows from the definition of the period. Resolving (\ref{Laplace}) row by row, we obtain $\pi_{i_j}=K-j$ for all vertices of the $(j+1)$-th row $(i_j,j), i_j=1,\dots,L,0\leq j \leq K-1$.

The structure of the $\pi$-Eulerian tour corresponding to the $(K+1)$-th period is illustrated in Fig.\ref{uni}.
The initial position $w$ of the chip is $S$. The chip-and-rotor state $(S,\rho_0)$ is the unicycle with the cycle $S,1,2,\dots,7,S^{'}$. The rings around the vertices denote multiple full rotations of the rotors during the tour. The number of rings around vertex $v$ is $\pi_v$. Each ring around vertex $v$ belonging to one of rows $2,3,\dots,K$ corresponds to four exits from $v$ along outgoing edges in a fixed, say, clockwise order. Each ring belonging to the first row corresponds to three exits to the left, up and to the right. The total number of directed edges in the graph $\overline{G_K }$ is $4K L$. The picture of rings for a given graph can be called a standard ring configuration. The determined components of the vector  $\boldsymbol{\pi}$ allow us to find the length of the $\pi$-Eulerian tour $\sum \pi_v d_v =L(2K^2+K+1)$ where summation is over all vertices of the graph $\overline{G_K }$ including vertex $S$.

If the chip-and-rotor state $(S,\rho_0)$ is a unicycle, then it is recurrent according to Theorem 3.8 in \cite{HLMPPW}. The orbit of a recurrent state is the set of of all states obtained from $(S,\rho_0)$ by iterating the rotor-router operation. It follows from above that the orbit starting from $(S,\rho_0)$ coincides with the set of states passed by the $\pi$-Eulerian tour starting from the same state. The number of orbits for an arbitrary strongly connected directed graph is found by Trung Van Pham in \cite{Pham} (see also \cite{FarrLev}).

\subsection{Multicycles}

The rotor configuration $\rho(v)$, $v\in G_K$ obtained after the $K$-th period is not necessarily loop-free and the chip-and-rotor state on the graph $\overline{G_K }$ is not necessarily a unicycle. The vertices of top row $(1,K-1),\dots,(L,K-1)$ of the aggregation cluster $A(N)$ are filled one by one with $L$ particles which remain immovable during the $K$-th period.
Therefore, the rotors at these vertices may have any of four orientations with equal probabilities. As a result, the rotor configuration $\rho(v)$, $v\in G_K$ can contain loops involving a part of vertices or all vertices of the top row.

In \cite{PPP15}, a concept of {\it multicycles} was proposed to deal with rotor configurations containing more than one loop. A multicycle $(w_1,w_2,\dots,w_k;\rho)$ is a chip-and-rotor state containing $k$ non-intersecting cycles $c_1,c_2,\dots,c_k$ and $k$ chips at vertices $w_1,w_2,\dots,w_k$ belonging to the cycles. The vertex $w_i$ is arbitrarily chosen from the vertices of the cycle $c_i, 1\leq i \leq k$.
A multicycle analog of the above mentioned Lemma 4.9 is the following. Let $G$ be an Eulerian directed graph with $m$ edges and $(w_1,w_2,\dots,w_k;\rho)$ be a multicycle in $G$. We iterate the rotor-router operation $m_1$ times starting from $w_1$ until the chip returns to $w_1$, the rotor $\rho(w_1)$ returns to its initial position and the rotor state is $\rho_1(v)$. Then, we iterate the rotor operation $m_2$ times starting from $w_2$ and from configuration $\rho_1(v)$ until the chip and the rotor $\rho_1(w_2)$ return to their initial positions in $\rho(v)$ and the rotor state is $\rho_2(v)$. Continuing with $w_3,w_4,\dots,w_k$, we iterate the initial state $m_1+m_2+\dots,m_k=m$ times and return to $(w_1,w_2,\dots,w_k;\rho)$. Each rotor makes one full turn and each edge is passed exactly once. The order of marking the vertices $w_1,w_2,\dots,w_k $ is inessential.

The properties of multicycles of Eulerian graphs can be readily generalized to non-Eulerian graphs by the same construction of multigraphs \cite{FarrLev} as in the case of unicycles. Having the non-Eulerian graph $\overline{G_K }$, we construct the auxiliary Eulerian graph $\widetilde{G_K}$ and apply the analog of the Lemma 4.9 for multicycles.
 \begin{figure}[!ht]
\includegraphics[width=90mm]{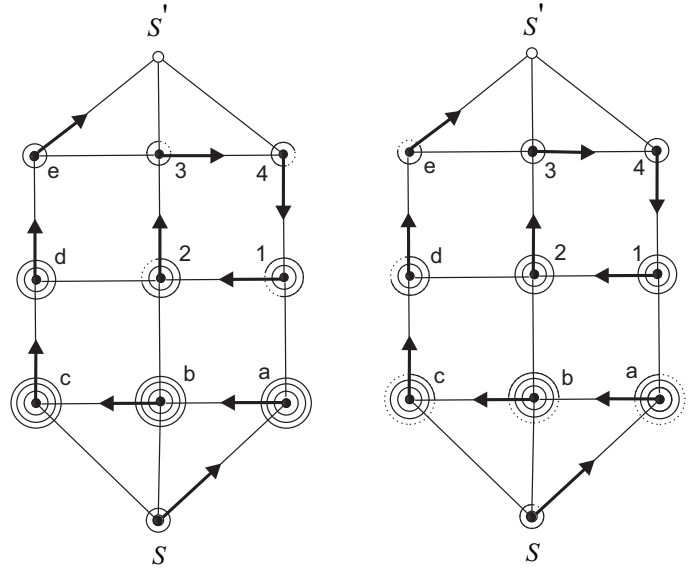}
\vspace{-1mm}
\caption{\label{multi}The graph  $\overline{G_K }$ for $K=3,L=3$.The chip-and-rotor state $(w_1,w_2;\rho_0)$ is a multicycle containing two cycles: $S,a,b,c,d,e,S^{'}$ and 1,2,3,4 with $w_1=S$ and $w_2=2$. The rings drawn by the full and dotted lines denote the rotor rotations during two $\pi$-Eulerian tours. The left side: The  first tour denoted by full lines is started from $w_1=S$, the second one denoted by dotted lines is started from $w_2=2$. The right side: The first tour is started from $w_2=2$, the second one from $w_1=S$.}
\end{figure}

An example of $\pi$-Eulerian
tours is shown in Fig.\ref{multi} where the graph $\overline{G_K }$ corresponds to the aggregation cluster $A(N)$ with $N=9$ obtained after 3 periods on the cylinder with $L=3$.
The initial chip-and-rotor state $(w_1,w_2;\rho_0)$ is a multicycle containing two cycles: $S,a,b,c,d,e,S^{'}$ and 1,2,3,4. Since the positions of chips on cycles are arbitrary, we chose for definiteness $w_1=S$ and $w_2=2$. The components of vector $\boldsymbol{\pi}$ are $\pi_S=1; \pi_e=\pi_3=\pi_4=1; \pi_d=\pi_2=\pi_1=2$ and $\pi_c=\pi_b=\pi_a=3$.
The first $\pi$-Eulerian tour starting in vertex $S$ and returning to $S$ after 62 steps is denoted in the left part of Fig.\ref{multi} by the full lines. The second $\pi$-Eulerian tour starting in vertex 2 and returning to 2 consists of 4 steps and is denoted by the dotted lines. Two tours together generate the standard ring configuration and contain $\sum_v \pi_vd_v=66$ steps.

One can notice that the rings around the vertices of the top row are not completed after the first tour. Nevertheless, the completed segments of these rings reach the edges directed to $S^{'}$ in accordance with the results of Section II. Indeed, the first tour coincides with the closed path on $\overline{G_K }$ traced by the chip starting from $S$ with the rotor configuration $\rho_0$. On the other hand, it is equivalent to $L$ paths traced by $L$ particles starting at $S$ during the $(K+1)$-th period. But according to Section II, all $L$
particles pass the edges directed upwards from the vertices of the top row exactly once.

The right part of Fig.\ref{multi} represents the case when the first $\pi$-Eulerian tour starts from the vertex 2. Again, the segments of the first tour are denoted by the full lines. The number of steps in the first tour is 58. The second tour started from $S$ is denoted by dotted lines and consists of 8 steps, which complete the first tour to the standard ring configuration.

For the general graph $\overline{G_K }$, the construction of multicycle $(w_1,w_2,\dots,w_k;\rho)$ allows one to determine an upper bound for the length of the single closed path  during the $(K+1)$-th period. This closed path is traced by the first of $k$ chips placed to $w_1=S$ in the initial chip-and-rotor state corresponding to the beginning of the $(K+1)$-th period. The other $(k-1)$ chips at $w_2,w_3,\dots, w_k$ generate auxiliary paths that show how far the first path is from completing the standard ring configuration. Each of these cycles takes at least two steps to return to its starting vertex $w_i, 2\leq i \leq k$ since the minimal cycle is a dimer formed by two opposite rotors. Taking into account that the total length of the $\pi$-Eulerian tours of all chips is $\sum_v \pi_vd_v$, we conclude that the length of the closed path starting at $S$ is bounded by $L(2K^2+K+1)-2(k-1)$ for the rotor configuration on $\overline{G_K }$ containing $k$ cycles.

\section{Discussion}

\begin{figure}[!ht]
\includegraphics[width=70mm]{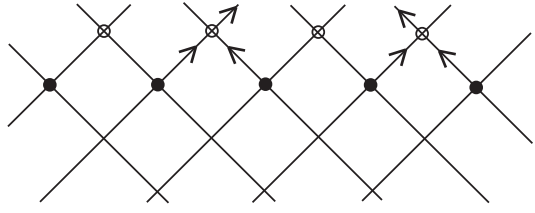}
\vspace{-1mm}
\caption{\label{turn}A part of cylinder obtained from a turned square lattice.
}
\end{figure}

Besides the rotor-router dynamics, the perfect surface obtained in Section II is a consequence of two factors: the cylinder geometry and the orientation of the square lattice. The crucial property of the considered aggregation process is the independent occupation by particles of a new row of vertices during each period. This property vanishes for lattices of different symmetry, for instance, for the honeycomb lattice or for the square lattices turning 45 degrees from the cylinder generatrix. In Fig.\ref{turn}, we show a part of cylinder taken out of a turned square lattice. As above, the aggregation process starts from the bottom of the cylinder.
Two particles from the row of vertices marked by  the filled circles can reach one vertex not visited before and marked by an open circle. One of the particles stops, according to the aggregation rules and the second continues a motion in a randomly chosen direction preventing the formation of a perfect surface. Then, the question about a constant width of boundary fluctuations remains open in this case.

At the same time, a proof of the constant width of the boundary for the square lattice having an arbitrary orientation on the cylinder would be an essential step towards the proof of the constant bonds for fluctuations of the circular form of the aggregation cluster on the two-dimensional square lattice. Indeed, the formation of the boundary of a large circle in the intervals corresponding to four cardinal directions resembles the aggregation process on the cylinder orientated along the lattice axes. A turned square lattice reproduces fluctuations of the circular form in an arbitrary direction when the cluster radius tends to infinity.
\section*{Acknowledgments}

The author acknowledges the support of the RFBR, grant 16-02-00252.


\begin{thebibliography}{99}

\bibitem{Meakin}P.Meakin and J.M.Deutch, The formation of surfaces by diffusion-limited annihilation, J.Chem.Phys. 85, 2320-2325 (1986).

\bibitem{Lawleretal} G.F.Lawler, M.Bramson and D.Griffeath, Inernal diffusion-limited aggregation, Ann.Prob. 20(4), 2117-2140 (1992).

\bibitem{Lawler} G.F.Lawler,Subdiffusive fluctuations for internal diffusion-limited aggregation, Ann.Prob. 23(1), 71-86 (1995).

\bibitem{AG1}A.Asselah and A.Gaudilli\`{e}re, From logarithmic to subdiffusive polynomial fluctuations for internal DLA and related
growth model, Ann.Prob. 41,1115-1159(2013).

\bibitem{AG2}A.Asselah and A.Gaudilli\`{e}re,Sub-logarithmic fluctuations for internal DLA, Ann.Prob. 41, 1160-1179 (2013).

\bibitem{JLS1}D.Jerison,L.Levine and S.Sheffield, Logarithmic fluctuations for internal DLA, J. Amer. Math. Soc. 25,271-301 (2012).

\bibitem{JLS2}D.Jerison,L.Levine and S.Sheffield,Internal DLA in higher dimensions, Electron.J.Prob.18, No 98 (2013).

\bibitem{Kleber}M.Kleber, Goldbug variations, Math.Intelligencer 27, 55-63 (2005).

\bibitem{PDDK} V.B. Priezzhev, D. Dhar, A. Dhar, and S. Krishnamurthy, Eulerian walkers as a model of self-organized
criticality, Phys. Rev. Lett. 77, 5079--5082 (1996).

\bibitem{HLMPPW}
A.E. Holroyd, L. Levine, K. Meszaros, Y. Peres, J. Propp and D.B. Wilson.
Chip-Firing and Rotor-Routing on Directed Graphs.
Progress in Probability 60, 331--364 (2008).

\bibitem{LP05}
L. Levine and Y. Peres. The rotor-router shape is spherical.
Math. Intelligencer 27(3), 9--11 (2005).

\bibitem{LP07}
L. Levine and Y. Peres.
Strong spherical asymptotics for rotor-router aggregation and the divisible sandpile.
Potential Analysis 30(1), 1--27 (2009).

\bibitem{LP16}
L. Levine and Y. Peres.
Laplacian growth, sandpiles and scaling limits.
arXiv:16611.00411v1 [mathPR].

\bibitem{JLS3}D.Jerison,L.Levine and S.Sheffield,Internal DLA for cylinders, Advances in Analysis:The Legacy of Elias M.Stein, chapter 8,
189-214, Prinston University Press (2014).

\bibitem{Kager} W.Kager, L.Levine, Rotor-router aggregation on the layered square lattice, Electron. J. Combinatorics 17, R152 (2010).

\bibitem{Pham} Trung Van Pham, Orbits of rotor-router operation and stationary distribution of random walks on directed graphs, Adv.Applied Math.70, 45--53, (2015).

\bibitem{FarrLev}M.Farrell and L.Levine, Multi-Eulerian tours of directed graphs, Electron.J.Combinatorics 23(2)P2.21 (2016)

\bibitem{PPP15}V.V.Papoyan, V S Poghosyan and V B Priezzhev, A loop reversibility and subdiffusion of the rotor--router walk, J. Phys. A: Math. Theor. 48, 285203 (2015).



\end{thebibliography}
\end{document}